\documentstyle[12pt,dina4]{article}
\newcommand{\ve}{{\varepsilon}}
\def\N{{{\rm I} \! {\rm N}}}
\def\R{{{\rm I} \! {\rm R}}}
\def\Z{\rm Z \hskip 2pt \llap Z}
\def\C{\!\rlap{\rm C}\hskip 2.9pt {\vrule height 6 pt width 0.7 pt}\,} 
\def\Q{\rm O\hskip 2pt \llap Q}
\newcommand{\K}{{\cal K}}
\newcommand{\T}{{\cal T}}
\newcommand{\B}{{\cal B}}
\newcommand{\cl}[1]{\overline{#1}}
\newcommand{\supp}{\rm supp}
\newtheorem{rem}{Remark}[section]
\newcommand{\remark}[1]{ 
\begin{rem}
#1
\end{rem}
}
\newtheorem{THEOREM}{Theorem}[section]
\newcommand{\thm}[1]{
\begin{THEOREM}
#1
\end{THEOREM}
}
\newtheorem{LEMMA}[THEOREM]{Lemma}
\newcommand{\la}[1]{
\begin{LEMMA}
#1
\end{LEMMA}
} 
\newtheorem{COROLLARY}[THEOREM]{Corollary}
\newcommand{\coro}[1]{
\begin{COROLLARY}
#1
\end{COROLLARY}
}
\newtheorem{PROPOSITION}[THEOREM]{Proposition}
\newcommand{\prop}[1]{
\begin{PROPOSITION}
#1
\end{PROPOSITION}
}
\newtheorem{DEFINITION}{Definition}[section]
\newcommand{\defi}[1]{
\begin{DEFINITION}
#1
\end{DEFINITION}
} 
\newcommand{\bull}{\rule{2mm}{2mm}}
\newcommand{\pf}{
\begin{flushleft}
{\bf Proof}
\end{flushleft}}
\newcommand{\nm}[1]{\parallel #1 \parallel}
\newcommand{\nmi}[1]{\nm{#1}_\infty}
\newcommand{\nmp}[1]{\nm{#1}_{\ell^p}}
\newcommand{\nmpi}[1]{\nm{#1}_{\ell^\infty}}
\newcommand{\nmc}[1]{\nm{#1}_{C[a,b]}}
\newcommand{\nmP}[1]{\nm{#1}_p} 
\newcommand{\nmQ}[1]{\nm{#1}_q}
\newcommand{\nmR}[1]{\nm{#1}_r} 
\newcommand{\nmu}[1]{\nm{#1}_{L^p_{u(\alpha )}}}
\newcommand{\nmw}[1]{\nm{#1}_{L^p_{v(\alpha )}}}
\newcommand{\Nmwd}[1]{\nm{#1}_{L^p_{w(2\alpha +2\delta +1)}}}
\newcommand{\nMu}[1]{\nm{#1}_{M^p_{u(\alpha )}}}
\newcommand{\nMw}[1]{\nm{#1}_{M^p_{w(\alpha )}}}
\newcommand{\nmU}[1]{\nm{#1}_{L^q_{u(\alpha )}}}
\newcommand{\nmW}[1]{\nm{#1}_{L^q_{w(\alpha )}}}
\newcommand{\nmUl}[1]{\nm{#1}_{L^q_{u(\alpha +2\lambda )}}}
\newcommand{\nmuk}[1]{\nm{#1}_{L^1_{u(\alpha )}}}
\newcommand{\nmwk}[1]{\nm{#1}_{L^1_{w(\alpha )}}}
\newcommand{\nmwg}[1]{\nm{#1}_{L^p_{v(\gamma )}}}
\newcommand{\nmgk}[1]{\nm{#1}_{L^1_{w(\gamma )}}}
\newcommand{\nMwg}[1]{\nm{#1}_{M^p_{\alpha ,\gamma }}}
\newcommand{\nmwr}[1]{\nm{#1}_{L^r_{w(\alpha )}}}
\newcommand{\nMwr}[1]{\nm{#1}_{M^{p,r}_{w(\alpha )}}}

\begin{document}
 \begin{center}
 {\LARGE{}On  weighted transplantation and multipliers for
 Laguerre expansions}\\[.4cm]
 {\sc{}Krzysztof Stempak \footnote{Institute of Mathematics, 
 University of Wroc\l aw, Wroc\l aw, Poland. The work of this author was 
 done during a stay at the Fachbereich Mathematik of the TH Darmstadt and was 
 supported in part by the TH Darmstadt and in part
 by the Deutsche Forschungsgemeinschaft under grant 436POL/115/1/0.
 .}
 and Walter  Trebels\footnote{Fachbereich Mathematik, TH Darmstadt, D--6100 
 Darmstadt, Germany.} }\\[.4cm]
  {(April 22, 1993 version)}\\[.4cm]
 \end{center}

 \bigskip
 {\bf Abstract.} 
Using the standard square--function method (based on the 
Poisson semigroup), multiplier conditions of H\"ormander type are derived 
for Laguerre expansions in $L^p$--spaces with power weights in the $A_p$-range; 
this result can be 
interpreted as an ``upper end point'' multiplier criterion which is fairly good 
for $p$ near $1$ or near $\infty $. A weighted generalization of 
Kanjin's \cite{kan} transplantation theorem 
allows to obtain a ``lower end point'' multiplier criterion 
whence by interpolation nearly ``optimal'' multiplier criteria (in dependance of 
$p$, the order of the Laguerre polynomial, the weight).


\bigskip
{\bf AMS(MOS) subject classifications.} 42A45, 42B25, 42C10

\bigskip 
\section{ Introduction}

In the last fifteen years a considerable activity initiated by a series of 
papers of Markett has taken place to study the 
Harmonic Analysis for Laguerre expansions. 

\smallskip \noindent
There are several ways of studying Laguerre expansions and corresponding 
multipliers. A principal one is based on the system $\{ L_k^\alpha \} 
_{k=0}^\infty $ of Laguerre polynomials, see Szeg\"o \cite[p. 100]{szego}, 
which are orthogonal in $L^2({\bf R}_+,\, x^\alpha e^{-x}dx),\; \alpha >-1.$ 
Weighted $L^p$-inequalities, proved for this system by Muckenhoupt \cite{mu}, 
turn out to be very useful and imply $L^p$-estimates for Laguerre function 
expansions. However, as far as multipliers are concerned, this system is hard 
to deal with. For instance, Askey and Hirschman proved that uniform boundedness 
of Ces\`aro multiplier sequences of any positive order is restricted to $p=2$ 
only.

\medskip \noindent
Another principal type of expansion was discussed by G\"orlich and 
Markett \cite{indag}, 
\cite{analy}. To become more precise, let us introduce the Lebesgue spaces
$$L^p_{v(\gamma )} = \{ f: \; \nmwg{f} = ( \int _0^\infty
|f(x) |^p x^{\gamma }\, dx )^{1/p} < \infty \} \; ,\quad 1 \le p <
\infty  ,\quad \gamma >-1.$$
If we define the  Laguerre function system $\{ l_k^\alpha \} $ by
$$l_k^\alpha (x) =( k!/\Gamma (k+\alpha +1) ) ^{1/2} e^{-x/2}L_k^\alpha (x),
\quad \alpha >-1,\quad n\in {\bf N}_0, $$
then it is orthonormal in $L^2({\bf R}_+,\, x^\alpha dx)$. If  
$  \gamma < p(\alpha +1)-1$ we can associate 
to $f\in L^p_{v(\gamma )}$  the Laguerre series
$$ f(x) \sim \sum _{k=0}^\infty a_k l_k^\alpha (x),\quad \quad 
a_k = \int _0^\infty f(x)l_k^\alpha (x)x^\alpha  dx. $$ 
Let $m = \{ m_k\} $ be a sequence of real or complex numbers and 
associate to $m$ the operator 
$$ T_m f(x) = \sum _{k=0}^\infty m_ka_k l_k^\alpha (x) $$
for all $f$ of the type $f(x)=p(x)e^{-x/2},\; p$ polynomial (the set of  these $f$ is 
dense in $L^p_{v(\gamma )}$, see \cite{mu}). The sequence $m$ is called 
a bounded multiplier on 
$L^p_{v(\gamma )}$,  notation $ m \in M^p_{\alpha ,\gamma }$, if
$$   \nMwg{ m } := \inf \{ C:\,  \nmwg{T_m f }\le C \nmwg{f} \} $$
is finite for all $f$ as before. We note the duality property $(1/p +1/p' 
=1)$ 
\begin{equation}\label{dual}
M^p_{\alpha ,\gamma } = M^{p'}_{\alpha , \alpha p'-\gamma p'/p}\, \, , \quad 
-1 < \gamma < p(\alpha +1) -1, \quad 1<p<\infty .
\end{equation}
The purpose of this paper is to obtain a better insight into the structure of 
the Laguerre multiplier space $M^p_{\alpha ,\gamma }$ from the point of view of 
sufficient conditions. We mention that Gasper and Trebels \cite{laguerre} 
discussed this question from the point of view of  necessary conditions. 

\bigskip \noindent
{\bf Remarks.} 1) In \cite{laguerre} another scaling of the orthogonal system 
$\{ l_k^\alpha \} $ is used, but it follows 
directly that our $M^p_{\alpha ,\gamma }$-space coincides with the space 
$M^p_{w(\gamma )}$ in \cite{laguerre}. The $ l_k^\alpha  $'s have the 
great advantage to possess a nice convolution structure as shown by G\"orlich 
and Markett \cite{indag}.

\medskip \noindent
2) The standard orthonormal system on $L^2({\bf R}_+,\, dx)$ is given by 
\begin{equation}\label{lcal}
 {\cal L}_k^\alpha (x) =( k!/\Gamma (k+\alpha +1)) ^{1/2} 
x^{\alpha /2}e^{-x/2}L_k^\alpha (x) ,\quad k\in {\bf N}_0 .
\end{equation}
Though the Parseval identity follows at once, its disadvantage is to be seen in 
the fact that no nice convolution structure is available. 
If we denote by ${\cal M}^p_{\alpha ,\gamma }$ the multiplier space with respect 
to the system $\{ {\cal L}_k^\alpha \}$ (introduced analogously to $M^p_{\alpha 
,\gamma }$), then one can easily show the following connection between the 
two multiplier spaces
\begin{equation}\label{mspaces}
{\cal M}^p_{\alpha ,\gamma } = M^p_{\alpha ,\gamma +\alpha p/2 } \; .
\end{equation} 

\medskip \noindent
3) Markett \cite{analy} and Thangavelu \cite{thangtrans} considered another 
orthonormalized system on $L^2({\bf R}_+,dx)$, namely $\{ \varphi _k^\alpha \}$,
$$ \varphi _k^\alpha (x) = {\cal L}_k^\alpha (x^2)(2x)^{1/2},\quad  
k\in {\bf N}_0 .$$ 
Again one makes sure that the associated multiplier space coincides with the 
previously introduced $M^p_{\alpha ,\gamma +\alpha p/2 +p/4-1/2}.$ 
Thus, if one wants to study multipliers with respect to these three orthonormal 
systems, one
only needs to discuss $M^p_{\alpha ,\gamma }$ for various parameters $\gamma $.

\bigskip \noindent 
The sufficient multiplier conditions we have in mind 
are described in terms of the following ``weak-bounded-variation'' sequence  
spaces from Gasper and Trebels \cite{wbv} 
$$ wbv_{q,s}= \{ m\in l^\infty : \, \| m \| _{q,s} < \infty \}, \quad 1 
\le q \le \infty, \quad s > 0, $$
where the norm on this sequence space is given by
$$ \| m \| _{q,s}= \| m \| _\infty +
\sup _n \left( \sum _{k=n}^{2n} k^{-1}|k^s\Delta^sm_k|^q\right) ^{1/q}, $$
(with the standard interpretation for $q=\infty $) 
and the difference operator $\Delta ^s$ of fractional order $s$ by
$$ \Delta ^s m_k = \sum _{j=0}^\infty A_j^{-s-1}m_{k+j}$$
whenever the sum converges. We are now ready to formulate a first 
sufficient multiplier criterion for Laguerre expansions (to 
be proved in Section 3).
\thm{ Let $\alpha \ge 0$, $1<p<\infty$,
and $m\in wbv_{2,s}$ for some $s>\alpha +1 $. Then 
$m\in M^p_{\alpha ,\gamma}$ and 
$\nMwg{ m } \le C \| m\| _{2,s}\, $, provided
\begin{equation}\label{wrestr}
 (\alpha +1) \max \{- p/2, \, -1\} < \gamma -\alpha <
    (\alpha +1) \min \{ p/2, \, p-1\}  ;
\end{equation}
i.e., under the above assumptions there holds in the sense of continuous 
embedding
$$ wbv_{2,s} \subset M^p_{\alpha , \gamma } \, . $$
}
{\bf Remarks.} 1) For $\gamma = \alpha \, $ Theorem 1.1 coincides with an 
unpublished result of Dietrich, G\"orlich, Hinsen, and Markett 
(due to a written communication of C. Markett). For integer $\alpha =n$, a 
weaker version of Theorem 1.1, namely 
$$ wbv_{\infty ,s(n)} \subset M^p_{n,n} ,\quad s(n) = \cases{ 
n+2 & odd $n$ \cr
n+3 & even $n$  } $$
can be read off from Thangavelu \cite{thanglec}.

\bigskip \noindent
2) By the weighted transplantation theorem to be deduced in Section 4 Theorem 
1.1 implies 
\begin{equation}\label{lep}
wbv_{2,s_0} \subset M^p_{\alpha ,\alpha }, \quad s_0>1, \quad \alpha \ge 0, \quad 
\frac{2\alpha +2}{\alpha +2} < p < \frac{2\alpha +2}{\alpha }\, .
\end{equation}
Following the lines of Connett and Schwartz \cite{cs},
interpolation between Theorem 1.1 for $\gamma = \alpha$ with $p$ near $1$ and 
(\ref{lep}) leads to Part a) of
\coro{
a) Let $\alpha \ge 0$ and $1<p<\infty$, then 
$$ wbv_{2,s} \subset M^p_{\alpha ,\alpha },\quad s>\max \{ s_c(p) ,\; 1 \} ,$$
where the quantity $s_c(p)=(2\alpha +2)|1/p -1/2|$ plays the role of a 
critical index in the  multiplier space $M^p_{\alpha ,\alpha }.$

\bigskip \noindent
b) Let $\alpha > -1/3$ and $1\le p < 2$, then 
$$ M^p_{\alpha ,\alpha } \subset wbv_{p',s'} , \quad 0<s'\le s_c(p) 
-4(1/p-1/2)/3. $$ 
} 

\bigskip \noindent 
Part b) is a slight extension of a result in Gasper and Trebels 
\cite{laguerre}. The corollary nicely shows where Laguerre multipliers live. 
Surprising is the smoothness gap of size $4|1/p-1/2|/3$, since in other 
settings a gap of size $|1/p-1/2|$ is well known, see e.g. \cite{jacobi} in the 
Jacobi case $(\alpha , -1/2)$.\\
If $s_c(p)\ge 1$ the result in Part a) is best possible in the following 
sense (continuous embedding)
$$ wbv_{2,s'} \not\subset M^p_{\alpha ,\alpha },\quad s'<s_c(p), $$
as the example of the Ces\`{a}ro means $m_{n,\nu}$ of order $\nu$ shows: 
observing 
$$ m_{n,\nu }(k) = \left\{ \begin{array}{l@{\quad,\, }l}
A_{n-k}^\nu /A_n^\nu & 0\le k\le n \\
0 & k>n
\end{array} \right.   , \quad 
\Delta^{s'}m_{n,\nu }(k) = \left\{ \begin{array}{l@{\quad,\, }l}
A_{n-k}^{\nu -s'} /A_n^\nu & 0\le k\le n \\
0 & k>n
\end{array} \right.  ,$$
it immediately follows that $m_{n,\nu } \in wbv_{2,s'}$ uniformly in $n$ if 
$s'<\nu +1/2 <s_c(p)$. But it is well known that
$\| m_{n,\nu } \| _{M^p_{\alpha ,\alpha }} \ge C (n+1)^{s_c(p)-\nu -1/2}$ (see 
e.g. \cite{cohen}). \\
On the other hand concerning the result in Part b) it is shown in 
\cite{laguerre} that at least for $p=1$ and $\alpha \ge 0$ the necessary 
condition is best possible. In \cite[II]{laguerre} there is given a modification 
of the $wbv$-condition (not nicely comparable in the $wbv$-framework) which 
suggests a smoothness gap of size $|1/p-1/2|$ if one could modify in the same 
sense the sufficient condition.

\bigskip \noindent
The plan of the paper is as follows. In Section 2 we develop the required 
square-function calculus based on the Poisson kernel. In Section 3 we give the 
relevant estimates for the Poisson means of the multiplier which allow 
one to use the 
standard square-function method (see e.g. \cite{cs}). In Section 4 we 
turn to a generalization of Kanjin's \cite{kan} transplantation theorem, thus 
including Thangavelu's \cite{thangtrans} modification, 
by admitting more general power weights. This allows some further insight 
into the structure of Laguerre multipliers. 

\bigskip \noindent
{\bf Acknowledgements.} Some of the ideas occurring in this paper were 
developed in other work done jointly by 
George Gasper and the second author. At this point we express 
our indebtedness to these discussions. The first author thanks the Fachbereich 
Mathematik of the  TH Darmstadt for hospitality and financial support.

\bigskip \noindent
\section{On square-functions for Laguerre expansions}
The main tool we will use is the twisted generalized convolution defined in 
$L^1({\bf R}_+,d\mu _\alpha ), \, d\mu _\alpha (x) = x^{2\alpha +1}dx,\; 
\alpha \ge 0, $ by 
$$ f\times g(x) = \int _0^\infty  \tau _xf(y)g(y)\, d\mu _\alpha (y),$$
where the twisted generalized translation operator $\tau _x$ is given by
$$\tau _x f(y) = \frac{\Gamma (\alpha +1)}{\pi ^{1/2}\Gamma (\alpha +1/2)}
\int _0^\pi f((x,y)_\theta){\cal J}_{\alpha -1/2}(xy\sin 
\theta) (\sin \theta )^{2\alpha } \, d\theta ,$$
${\cal J}_{\beta }(x) = \Gamma (\beta +1) J_\beta (x)/(x/2)^\beta $, 
$J_\beta $ denoting the Bessel function of order $\beta >-1$, and 
$$ (x,y)_\theta = (x^2 +y^2 -2xy \cos \theta)^{1/2};$$
this convolution is commutative -- for all this 
see G\"orlich and Markett \cite{indag} and Stempak \cite{stcon}.
For the following it is convenient to work with the transformed system
$$\psi _k^\alpha (x) = (2k!/\Gamma (k+\alpha +1))^{1/2} e^{-x^2/2}L_k^\alpha 
(x^2), \quad k\in {\bf N}_0, $$
which, on account of the orthonormality of $\{ l_k^\alpha \}$, is obviously 
orthonormal on $L^2({\bf R_+}, d\mu _\alpha )$. We introduce the 
norm 
$$\| f\| _{p,\delta } = \left( \int _0^\infty |f(x)|^px^{2\delta} d\mu _\alpha 
(x) \right) ^{1/p}; $$
this will not lead to any confusion with the $wbv$-norm.
We note that the $\psi _k^\alpha $'s are eigenfunctions, with eigenvalues 
$\lambda _k$, of the positive symmetric in $L^2 ({\bf R_+}, d\mu _\alpha )$  
differential operator $L$,
$$ L=-\left( \frac{d^2}{dx^2} +\frac{2\alpha +1}{x} \frac{d}{dx} -x^2 \right) ,
\quad \quad \lambda _k = 4k+2\alpha +2.$$
We mention the following transform property of the twisted generalized 
convolution with respect to $\psi _k^\alpha $-Laguerre expansions; 
if $f\sim \sum c_k\psi _k^\alpha $ and 
$f\times g \sim \sum c_kd_k\psi _k^\alpha $, then 
$g(x)\sim \Gamma (\alpha +1) \sum d_kL_k^\alpha (y^2)e^{-y^2/2}$. 
To $f\sim \sum c_k\psi _k^\alpha $ we associate its Poisson means 
$$P^tf = \sum _{k=0}^\infty \exp (-t\lambda _k)c_k \psi _k^\alpha =f\times p_t , 
\quad t>0, $$
with Poisson kernel $p_t(y)= c_\alpha \sum \exp (-t\lambda _k) L_k^\alpha 
(y^2)e^{-y^2/2}$, and the square-function
\begin{equation}\label{g1}
g_1(f)^2(x) = \int _0^\infty | \frac{\partial }{\partial t} P^tf(x) |^2t\, dt =
\int _0^\infty | \sum_{k=0}^\infty \lambda _k \exp (-t\lambda 
_k)c_k \psi _k^\alpha (x)| ^2 t\, dt
\end{equation}
for appropriate $f$.
Since for $\alpha \ge 0$ the semigroup $\{ P^t\} _{t>0} $ forms a positive 
contraction 
semigroup (see \cite{indag} and \cite{stpos}), by the Coifman, Rochberg, and 
Weiss refinement of  Stein's general Littlewood--Paley theory 
(see Meda \cite{meda}) 
\begin{equation}\label{basicg}
 C^{-1} \| f \| _{p,0}\le \| g_1(f)\| _{p,0}  \le C 
\| f \| _{p,0} , \quad 1<p<\infty 
\end{equation}
is true. We want to extend (\ref{basicg}) to the weighted case $\delta \neq 0$.
\prop{
Let $1<p<\infty ,\; \alpha \ge 0$, and $-(\alpha +1) < \delta <(\alpha 
+1)(p-1). $ Then
$$ C^{-1} \| f \| _{p,\delta }\le \| g_1(f)\| _{p,\delta }  \le C 
\| f \| _{p,\delta }.$$
}
{\bf Proof.}
Assume for the moment that the right inequality holds for the $\delta 
$'s indicated. We show that then the left one follows. 
Since by straightforward calculation, based on Parseval's identity, 
$\| f\| _{2,0} = C(\alpha ) \| g_1(f)\| _{2,0} \, $,
polarization and  H\"older's inequality give 
\begin{eqnarray*}
\left| \int _0 ^\infty f_1(x)f_2(x)d\mu _\alpha (x) \right| &\le  & C 
\int _0 ^\infty g_1(f_1)(x)x^{2\delta /p} x^{-2\delta /p} 
g_1(f_2)(x)d\mu _\alpha (x)\\
 & \le & C \| g_1(f_1)\| _{p,\delta } \| g_1(f_2)\| _{p',-\delta p'/p} .
\end{eqnarray*}
Setting $f_2(x) = x^{2\delta /p} h(x)$ and using the right hand inequality 
in Proposition 2.1 one obtains 
$$\left| \int _0 ^\infty f_1(x) x^{2\delta /p} h(x) d\mu _\alpha (x) 
\right| \le C \| g_1(f_1)\| _{p,\delta } \| h\| _{p',0} $$
for all $h\in L^{p'}(d\mu _\alpha )$ and 
$-(\alpha +1) < -2\delta p'/p <(\alpha +1)(p'-1) $ so that finally the 
converse of H\"older's inequality gives the desired  left hand side inequality 
$$\| f_1 \| _{p,\delta } =\| f_1 x^{2\delta /p} \| _{p,0}  \le C
\| g_1(f_1)\| _{p,\delta }.  $$

\medskip \noindent
To extend the right hand side inequality of (\ref{basicg}) to the weighted case 
$\delta \neq 0$ we adapt an approach of Stein \cite{sweight}.
The following is a variation of a corresponding remark in \cite[p. 271]{sbook}.
\la{
For $\alpha \ge 0$ 
let  $K(x,y)$ be  a homogeneous kernel, 
$K(\lambda x,\lambda y)=\lambda ^{-(2\alpha +2)}K(x,y)$, 
satisfying
$$ \int _0^\infty |K(1,y)|y^{-(2\alpha +2)/p}d\mu _\alpha(y) < \infty . $$
Then the operator
$$Tf(x)=\int _0^\infty K(x,y)f(y)\, d\mu_\alpha (y) $$ 
is bounded on $L^p(d\mu _\alpha ) $.
}
The particular kernel in the next lemma is the key of the desired extension.
\la{
For $\alpha \ge 0$ the kernel 
$$ K(x,y)= | 1-(x/y)^{2\delta /p} | \, 
\int _0^\pi (x,y)_\theta ^{-(2\alpha +2)} (\sin \theta )^{2\alpha } \, d\theta 
$$
satisfies the properties of Lemma 2.2 provided
$-(\alpha +1) < \delta < (\alpha +1)(p-1).$ 
}

{\bf Proof.} The homogeneity property of the kernel is clear. For the second 
property note that there are three singularities: $0,\, \infty $,  and $1$ and 
the  required integrability of $K(1,y)$ follows once we show that
\begin{equation}\label{a}
\int _0^\pi (1,y)_\theta ^{-(2\alpha +2)} (\sin \theta )^{2\alpha } \, d\theta 
\le C 
\cases{
1 & as  $y\to 0 $ \cr
|1-y|^{-1} & as $y\to 1$\cr
y^{-(2\alpha +2)} & as $y\to \infty .$
} 
\end{equation}
For small $y$ note that $(1,y)_\theta \approx 1,\, 0 < \theta <\pi $
whereas $(1,y)_\theta \approx y$ for large  $y$. Finally, 
for $y\approx 1$ it follows that
$$\int _0^\pi \frac{(\sin \theta )^{2\alpha }}{((1-y)^2+4y\sin 
^2(\theta/2))^{\alpha +1}} \, d\theta
\le C |1-y|^{-(2\alpha +2)}\int _0^{|1-y|} (\sin \theta )^{2\alpha }\, d\theta 
$$
$$\quad \quad \quad \quad +C \int _{|1-y|}^\pi \sin ^{-2} (\theta /2)\, d\theta 
\le C |1-y|^{-1} $$
so that (\ref{a}) is obvious.

\bigskip \noindent
The modification to twisted generalized convolution operators
of Stein's \cite{sweight} result now reads

\la{
Suppose $|A(x)| \le C x^{-(2\alpha +2)}$ and $Tf = A\times f$ is bounded on 
$L^p(d\mu _\alpha ),\; 1<p<\infty $. Then we also have the weighted inequality
$$ \| Tf\| _{p,\delta } \le C \| f\| _{p,\delta }, \quad 
-(\alpha +1) < \delta < (\alpha +1)(p-1).$$ 
}

\bigskip \noindent 
The proof follows along the pattern of \cite{sweight}. We only note that
$$ x^{2\delta /p} A\times f(x)= A\times (y^{2\delta /p}f)(x) -
\int _0^\infty (1-(x/y)^{2\delta /p})\tau _xA(y)f(y) y^{2\delta /p}\, 
d\mu _\alpha (y)$$
and, since  $|{\cal J}_{\alpha -1/2}(t)| \le 1, \; \alpha \ge 0$, 
the last integral can be estimated by 
$$ \int _0^\infty K(x,y)|f(y)| y^{2\delta /p}\, d\mu _\alpha (y),$$
where $K$ is the kernel from Lemma 2.3, hence the assertion by Lemma 2.2.

\bigskip \noindent
The above scalar-valued result can be extended to the case of functions taking 
their values in a Hilbert space, see \cite[pp. 45]{sbook}, by a repetition of 
the arguments given for the scalar-valued case.
\prop{
Let $A(x)$ be a function on ${\bf R}_+$ taking values in $B({\cal H}_1,
{\cal H}_2)$, ${\cal H}_i$ being separable Hilbert spaces, and satisfying 
$\| A\| \le C x^{-(2\alpha +2)}$. Further, if one defines $Tf= A\times f$
for $f\in L^p({\bf R}_+,{\cal H}_1) $ and if 
$$ \int _0^\infty \| Tf(x)\| ^p_{{\cal H}_2}x^{2\delta}d\mu _\alpha (x) 
\le C \int _0^\infty \| f(x)\| ^p_{{\cal H}_1}x^{2\delta}d\mu _\alpha (x) ,\quad 
1<p<\infty ,$$
holds for $\delta =0$, then the same inequality is valid for all $\delta $'s
satisfying $-(\alpha +1) < \delta < (\alpha +1)(p-1).$ 
}
To complete the proof of Proposition 2.1 we choose in the preceding proposition 
${\cal H}_1  = {\bf C}, \; {\cal H}_2 = L^2({\bf R}_+, t\, dt)$ and 
$Tf = (\partial /\partial t) P^t f= f\times (\partial /\partial t) p_t $. 
It is shown in Thangavelu \cite[Lemma 3.1]{thanglec} that 
$$ \| (\partial /\partial t) p_t (x) \| 
_{L^2({\bf R}_+, t\, dt)} \le C x^{-(2\alpha +2)} ,$$
where $ (\partial /\partial t) p_t =c_\alpha \sum  \lambda _k \exp 
(-t\lambda _k) L_k^\alpha (y^2)e^{-y^2/2}$.
Since (\ref{basicg}) holds, all hypotheses of 
Proposition 2.5 are satisfied and hence Proposition 2.1 is established.

\bigskip \noindent
For the proof of the multiplier Theorem 1.1 we need the standard variations on 
the $g_1$-function, namely $g_\sigma $-functions and a $g_\lambda 
^*$-function.\\
We note that a substitution in (\ref{g1}) gives
$$ g_1(f)^2(x) = \int _0^1| \sum_{k=0}^\infty \lambda _k r^{\lambda _k}c_k 
\psi _k^\alpha (x)| ^2 |\log r| \, \frac{dr}{r} . $$
Denoting $u(x,r) =P^tf(x),\; e^{-t} =r,$ and following Strichartz 
\cite{strich} we introduce as $\sigma $-th derivative of $u(x,r)$
$$     d_\sigma u(x,r) =  \sum_{k=0}^\infty \lambda _k^\sigma r^{\lambda _k}c_k 
\psi _k^\alpha (x)$$
and set 
\begin{equation}\label{gs}
 g_\sigma(f)^2(x) = \int _0^1|d_\sigma u(x,r) | ^2 |\log r|^{2\sigma -1}
 \, \frac{dr}{r} . 
\end{equation}
For the definiton of the $g_\lambda ^*$-function we need a generalized 
Euclidean translation (occurring in the framework of modified Hankel 
transforms, cf. \cite{stcon})
$$\tau ^E_x f(y) = \frac{\Gamma (\alpha +1)}{\pi ^{1/2}\Gamma (\alpha +1/2)}
\int _0^\pi f((x,y)_\theta)(\sin \theta )^{2\alpha } \, d\theta $$
and its associated convolution 
$$f*g(x) = \int _0^\infty \tau ^E_x f(y)g(y) \, d\mu _\alpha (y);$$
then the $g_\lambda ^*$-function is defined by 
$$ g_\lambda ^*(f)^2(x) = \int _0^1 K_{|\log r |}*|d_1u(\cdot ,r) | ^2
(x) |\log r|\, \frac{dr}{r} , $$
where $ K_t(y)=\delta _{\sqrt{t}}K(y),\; K(y) = (1+y^2)^{-\lambda }$ and 
$\delta _uf(y)= u^{-2(\alpha +1)}f(y/u)$ is an $L^1(d\mu _\alpha )$-invariant 
dilation.
\prop{   
a) $$ g_\rho (f)(x) \le C g_\sigma (f)(x) \quad a.e., 1\le 
\rho \le \sigma $$
for all $f\in L^p(x^{2\delta }d\mu _\alpha )$ for which the right hand side 
makes sense.

\medskip \noindent
b) $$ \| g_\lambda ^*(f)\| _{p, \delta } \le C \| f\| _{p,\delta },\quad 
\lambda > \alpha +1 ,\quad p \ge  2$$ 
provided  $-(\alpha +1)< \delta < p(\alpha +1)(1/2-1/p).$
}

{\bf Proof.} Assertion a) is proved just as in Strichartz \cite{strich}. For 
the proof of b) we note that the method in Stein \cite[p. 91]{sbook} works. 
Start with the basic inequality ($M$ denotes the Hardy-Littlewood maximal 
operator in the homogeneous space $({\bf R}_+,d\mu _\alpha , \rho )$ where 
$\rho $ is the usual distance on ${\bf R}_+$).
\begin{equation}\label{basic} 
\int _0^\infty g_\lambda ^*(f)^2(x)h(x)d\mu _\alpha \le C 
\int _0^\infty g_1(f)^2(x) Mh(x)d\mu _\alpha ,
\end{equation}
 which also holds in the present setting on account of the formulae 
(3.6), (3.8) and (3.13) in Stempak \cite{stcon}. Here the assumption $\lambda > 
\alpha +1$ implies $K(y)\in L^1(d\mu _\alpha )$.

\medskip \noindent
In the case $p=2$ choose $h(x)= x^{2\delta }$ in (\ref{basic}) and note that 
the maximal function 
$Mh(x) = C x^{2\delta }Mh(1)$ for $-(\alpha +1) <\delta \le 0$, hence the 
assertion by Proposition 2.1.\\
In the case $q=p/2>1$ choose $h(x)=x^{4\delta /p}h_1(x)$ in (\ref{basic}) 
and apply H\"older's inequality ($ 1/q+1/q'=1$) to obtain
$$\int _0^\infty g_\lambda ^*(f)^2(x) x^{4\delta /p} h_1(x)d\mu _\alpha \le C 
\int _0^\infty g_1(f)^2(x) x^{4\delta /p}x^{-4\delta /p}  M(h_1x^{4\delta /p})
d\mu _\alpha $$
$$ \le C \left( \int _0^\infty g_1(f)^p(x)x^{2\delta} d\mu _\alpha  \right) 
^{2/p} \left( \int _0^\infty M(h_1x^{4\delta /p})^{q'}(x)x^{-4q'\delta /p}
d\mu _\alpha \right) ^{1/q'}$$
$$\le C \| g_1(f)\| _{p,\delta }^{2} \| h_1\| _{q',0}
\quad \quad \quad \quad \quad \quad \quad \quad \quad \quad \quad \quad ,$$
for $x^{-4q'\delta /p} \in A_{q'}(d\mu _\alpha )$  by the assumption on $\delta $  (see 
\cite[II]{stcon}). Taking the supremum over all $h_1$, 
$\| h_1\| _{q',0} \le 1$ gives 
$$\| g_\lambda ^*(f)\| _{p,\delta }^2 
= \| g_\lambda ^*(f) x^{4\delta /p} \| _{q,0}
 \le C \|   g_1(f)\| _{p,\delta } ^{2} \le C \| f\| _{p,\delta } ^{2} .$$

\section{Proof of Theorem 1.1.}
Since we follow the standard method (see e.g. \cite[p.73]{cs}), we only 
indicate the main steps. We use the notation
$$f\sim \sum c_k\psi _k^\alpha ,\; \quad S_mf\sim \sum m_kc_k\psi _k^\alpha ,$$ 
and work on a dense subset of $L^p(x^{2\delta }d\mu _\alpha )$ 
(see \cite{mu}) such that we 
can write $=$ instead of $\sim $ in the preceding formulae. We note that the 
multiplier space associated to $\| \cdot \| _{p,\delta }$ coincides with 
$M^p_{\alpha ,\alpha +\delta }$.  Thus all we need to show, under the 
assumptions of Theorem 1.1, is 
\begin{equation}\label{mesti}
\| S_mf\| _{p,\delta } \le C \| m \| _{2,s} \| f\| _{p,\delta } , \quad \delta 
=\gamma - \alpha .
\end{equation}
If we assume for the moment that
\begin{equation}\label{pest}
g_{s+1}(S_mf)(x) \le C \| m \| _{2,s} g_s^*(f)(x) \quad a.e.
\end{equation}
holds, then (\ref{mesti}) is proved in the case $p\ge 2, \; 
-(\alpha +1)< \delta < p(\alpha +1)(1/2-1/p)$ by the 
following chain of norm inequalities if we choose 
 $\lambda = s > \alpha +1$ in Proposition 2.6 b)
$$ \| S_mf\| _{p,\delta } \le C \| g_1(S_mf)\| _{p,\delta } \le C 
\| g_{s+1}(S_mf)\| _{p,\delta } $$
$$\le C \| m \| _{2,s} \| g_s^*(f)\| _{p,\delta } \le C \| m \| _{2,s} 
\| f\| _{p,\delta } .$$
In the case $p=2$  the result (\ref{mesti}) extends at once to 
$-(\alpha +1)< \delta < \alpha +1$ by duality (\ref{dual}). Repeating 
the interpolation and duality arguments in Hirschman \cite[p.50]{hirsch} 
yields Theorem 1.1 provided (\ref{pest}) holds.

\bigskip \noindent
Let us turn to the proof of (\ref{pest}). First we note that (up to a constant) 
we don't change the $g_\sigma $-function (\ref{gs}) if we substitute $r^2$ for
$r$. By the properties of twisted convolution we have 
$$d_{s+1}S_mu(x,r^2)= d_sM(\cdot ,r) \times d_1u(\cdot ,r)(x) =
\int _0^\infty \tau _x d_1u(y,r) d_sM(y,r)d\mu _\alpha (y),
$$ 
where
$$ d_sM(y,r) =  c_\alpha \sum _{k=0}^\infty \lambda _k ^s m_k r^{\lambda _k} 
L_k^\alpha (y^2)e^{-y^2/2} .$$
Basic properties of $ d_sM(y,r)$, we need for the proof of (\ref{pest}), 
are contained in 
\prop{
Let $\alpha \ge 0$ and $m\in wbv_{2,s}$ for $s > \alpha +1$. Then
\begin{itemize}
\item[a)] $$\sup _y|d_sM(y,r)| \le C r^{2\alpha +2} (1-r)^{-s-\alpha -1}
\| m\| _\infty ,$$
\item[b)] $$ \int _0^\infty |y^s d_sM(y,r)|^2d\mu _\alpha (y) \le C 
r^{4\alpha +4} (1-r)^{-s-\alpha -1} \| m\| ^2_{2,s}. $$
\end{itemize}
}
Suppose that Proposition 3.1 is proved, then one obtains 
by the Cauchy--Schwarz inequality 
$$g_{s+1}(S_mf)^2(x)= \int _0^1
|\left( \int _0^{\sqrt{1-r}} + \int _{\sqrt{1-r}}^\infty \right) 
\tau _x d_1u(y,r) d_sM(y,r)d\mu _\alpha (y)
|^2 |\log r|^{2s+1}\frac{dr}{r}$$
$$ \le C \| m\| ^2_\infty \int _0^1 \frac{r^{4\alpha +4}}{(1-r)^{\alpha +2s+1}}
\int _0^{\sqrt{1-r}} |\tau _x d_1u(y,r)|^2 d\mu _\alpha (y) 
|\log r|^{2s+1}\frac{dr}{r}$$
$$\quad \quad \quad 
+ C \| m\| _{2,s}^2 \int _0^1 \frac{r^{4\alpha +4}}{(1-r)^{\alpha +s+1}}
\int _{\sqrt{1-r}}^\infty |y^{-s}
\tau _x d_1u(y,r)|^2 d\mu _\alpha (y) |\log r|^{2s+1}\frac{dr}{r}$$
$$ \le C \| m\| ^2_\infty \int _0^1 \frac{r^{4\alpha +4}|\log r|^{2s}}
{(1-r)^{\alpha +2s+1}}
\int _0^{\sqrt{1-r}} \tau ^E_x (|d_1u(y,r)|^2) d\mu _\alpha (y) 
|\log r|\frac{dr}{r}$$
$$\quad \quad \quad 
+ C \| m\| _{2,s}^2 \int _0^1 \frac{r^{4\alpha +4}|\log r|^{2s}}
{(1-r)^{\alpha +s+1}}
\int _{\sqrt{1-r}}^\infty y^{-2s}\tau ^E_x (|d_1u(y,r)|^2) d\mu _\alpha (y) 
|\log r|\frac{dr}{r}$$
because $|\tau _x f(y)| \le \tau ^E_x(|f|)(y)$, see \cite{stcon}. Now we use
in the first integral $y^2/|\log r| \le 1 $ if $0<y <\sqrt{1-r}$ and 
$r^4 (|\log r|/(1-r))^{2s+\alpha +1} \le  C$ for $0<r<1$, and in the second 
$r|\log r|^s y^{-2s} \le C (1+ y^2/|\log r|)^{-s}$ if $1-r < y^2$   
and arrive at 
\begin{eqnarray*}
g_{s+1}(S_mf)^2(x) 
 & \le & C \{ \| m\| ^2_\infty +\| m\| _{2,s}^2 \} 
\int _0^1 K_{|\log r|} * |d_1u(\cdot ,r)|^2(x) |\log r|\frac{dr}{r} \\
 & \le & C \| m\| _{2,s}^2 g_s^*(f)^2(x)
\end{eqnarray*}
by the definition of the $wbv$-norm.

\bigskip \noindent
Thus there remains only to prove Proposition 3.1. On account of 
\cite[Lemma 1]{analy}, 5th case, there holds
$$\sup _y |L_k^\alpha (y^2)e^{-y^2/2}| \le C (k+1)^\alpha $$
and hence Part a) 
$$\sup _y|d_sM(y,r)| \le C \| m\| _\infty r^{2\alpha +2}\sum _{k=0}^\infty 
(k+1)^{s+\alpha }r^{4k}\le C \| m\| _\infty r^{2\alpha +2} (1-r)^{-s-\alpha 
-1},$$
since $\sum (k+1)^\gamma r^k \le C (1-r)^{-\gamma -1}$ is true 
for $\gamma >-1 $.\\
To prove Part b) we use a weighted Parseval formula. First we note that the 
coefficients $c_k$ in the $\psi _k^\alpha $-expansion are related to the 
Fourier Laguerre coefficients ${\hat g}_\alpha (k)$ of \cite{laguerre} 
in the following way 
$$c_k = \Gamma (\alpha +1)(\Gamma (k+\alpha +1)/2k!)^{1/2} 
\left( e^{y/2}f(y^{1/2})\right) _\alpha  {\hat{\ }(k)}, \quad k\in {\bf N}_0;
$$
then it  follows from Gasper and Trebels 
\cite[I]{laguerre}, formulae (3) and (5) there, that
$$\sum _{k=0}^\infty A_k^{\alpha +s}|\Delta ^s ( c_k \sqrt{k!}/
\sqrt{\Gamma (k+\alpha +1)})|^2
\approx \int _0^\infty |f(x)|^2x^{2\alpha +2s+1} dx,$$
hence for the particular case $f=d_sM(\cdot ,r)$
$$ \int _0^\infty |y^sd_sM(y,r)|^2 d\mu_\alpha (y) 
\le C \sum _{k=0}^\infty A_k^{\alpha +s} |\Delta ^s (\lambda _k ^sm_kr^{\lambda 
_k} )|^2 =:I.$$
We have to dominate $I$.
Since similar computations for integer $s$ are contained in \cite[p. 69]{cs}
we only sketch the proof in that case. First note that
$\Delta ^\kappa r^k = (1-r)^\kappa r^k, \; \kappa >0$, and that 
$$ |\Delta ^j (\lambda _k ^sr^{4k} )|^2 \le C \sum _{i=0}^j (1-r)^{2i} r^{8k} 
(k+1)^ {2(s-i-j)};$$ 
then use these formulae in Leibniz' formula for differences 
$$I \le C r^{4\alpha +4} \sum _{j=0}^s \sum _{k=0}^\infty
A_k^{\alpha +s} |\Delta ^jm_k|^2|\Delta ^{s-j} (\lambda _{k+j}^sr^{4k+4j})|^2
\quad \quad \quad $$
$$ \le C r^{4\alpha +4} \sum _{j=0}^s \| m\| ^2_{2,j}(1-r)^{-\alpha -s-1}
\le C r^{4\alpha +4}\| m\| ^2_{2,s}(1-r)^{-\alpha -s-1}$$
by the embedding properties of the $wbv$-spaces, see \cite{wbv}.\\
If $s$ is strictly fractional, similar computations have been carried through 
in the proof of \cite[Lemma 1]{stud} and again we only sketch the proof. 
We use Peyerimhoff's \cite{peyer} version of Leibniz' 
formula for fractional differences which in our instance reads:\\
Let $s= [s]+\kappa ,\; 0< \kappa < 1$. Then
$$ \Delta ^s (\lambda _k ^sr^{\lambda _k} m_k)  
  =  r^{2\alpha +2}\sum _{i=0}^{[s]} {s \choose i}
\Delta ^i (\lambda _k ^sr^{4k})\Delta ^{s-i} m_{k+i} + r^{2\alpha +2}
m_k\Delta ^s (\lambda _k ^sr^{4k})$$
$$ + (-1)^{[s]}r^{2\alpha +2}R_k,$$
where the remainder term $R_k$ is given by
$$  \sum _{i=k+1+[s]}^\infty
A_{i-k}^{-s-1}(m_i-m_k)\sum _{j=k+1}^{i-[s]} A_{i-[s]-j}^{-[s]-1} \{ 
\Delta ^{[s]}(\lambda _j^sr^{4j}) - \Delta ^{[s]}(\lambda _k ^sr^{4k})\} . $$
Up to the terms which contain $\Delta ^{s-i} m_{k+i}$ no smoothness of the 
sequence $m$ is required and in that case $|m_k|$ can be crudely estimated by
$\| m\| _\infty $. To give an idea of the type of analysis required let us look 
at
$$|\Delta ^s (\lambda _k ^sr^{4k})| \le C \Delta ^{s-[s+1]}\sum _{l=0}^{[s+1]}
|\Delta ^lr^{4k}| \, |\Delta ^{[s+1]-l} \lambda _k ^s | \quad \quad \quad 
\quad \quad \quad $$
$$\le C \sum _{l=0}^{[s+1]}(1-r)^l \left( \sum _{i=0}^k + \sum _{i=k+1}^\infty 
\right) A_i^{[s]-s}r^{4(k+i)}(k+i+1)^{s-[s+1]+l} .$$
For $0 \le i \le k$ one has $(k+i+1) \approx (k+1)$ and thus 
$\sum _0^k \ldots \le C r^{4k}(k+1)^{l-1}$; for $i >k$ one can replace 
$(k+i+1)$ by $(i+1)$ and one obtains
$$\sum _{i=k+1}^\infty A_i^{[s]-s}r^{4(k+i)}(k+i+1)^{s-[s+1]+l}
\le C r^{4k} \left\{ \begin{array}{l@{\, ,\, }l}
r^{4(k+1)}((k+1)(1-r^{4k}))^{-1/2} & l=0\\
r^{4(k+1)}(1-r^{4k})^{-l} & l\neq 0, 
\end{array} \right.
$$
hence
$$\sum _{k=0}^\infty A_k^{\alpha +s}|r^{2\alpha +2}m_k\Delta ^s (\lambda _k 
^sr^{4k})|^2 \le C r^{4\alpha +4}(1-r)^{-\alpha -s-1}\| m\| _{2,s}^2. $$
One should note that in these estimates essentially 
$s>\alpha +1 \ge 1/2 $ is used. The contribution of the error term is estimated 
analogously.\\
Let us conclude with estimating the terms which require smoothness of the 
multiplier sequence in question. Observe that $(k+1+i) \approx (k+1)$ if  $0 
\le i < s$. Then
$$ \sum _{k=1}^\infty |r^{2\alpha +2} \Delta ^i (\lambda _k ^s r^{4k})
\Delta ^{s-i} m_{k+i}|^2 \quad \quad \quad \quad \quad \quad \quad \quad \quad 
\quad \quad \quad $$
$$\le C \sum _{n=1}^\infty \sum _{k=2^{n-1}}^{2^n-1} (k+1)^{\alpha +s}
\sum _{l=0}^i (1-r)^{2l}r^{8k}(k+1)^{2l+1} |(k+i)^{s-i}\Delta ^{s-i}m_{k+i}|^2
(k+i)^{-1} $$
$$\le C \sum _{l=0}^i (1-r)^{2l} \| m\| _{2,s-i}^2\sum _{n=1}^\infty 
2^{n(\alpha +s+2l+1)} r^{4\, 2^n} \quad \quad \quad \quad \quad \quad \quad $$
$$ \le C \| m\| _{2,s-i}^2(1-r)^{-\alpha -s-1} \le C 
\| m\| _{2,s}^2(1-r)^{-\alpha -s-1} \quad \quad \quad \quad \quad $$
again by the embedding properties of the $wbv$-spaces. \\
Thus Proposition 3.1 holds and Theorem 1.1 is established.

\section{A weighted transplantation theorem.}
In this context it is convenient to work with the 
Laguerre functions $\{ {\cal L}_k^\alpha \} $, introduced in (\ref{lcal}) at 
the beginning.
The following transplantation theorem for Laguerre function expansions has 
recently been proved by Kanjin \cite{kan}.

\bigskip \noindent
{\bf Theorem A.} {\it Let $\alpha ,\, \beta >-1$ and $\varepsilon = \min \{
\alpha ,\beta \}.$ If $\varepsilon \ge 0$, then 
$$\| \sum b_k {\cal L}_k^\alpha \| _{L^p({\bf R}_+,dx)}
\le  C \| \sum b_k {\cal L}_k^\beta \| _{L^p({\bf R}_+,dx)} $$
for $1<p<\infty $, where $C $ is a constant independent of $f$. If 
$-1 < \varepsilon <0$, then the assertion remains true provided $p$ satisfies 
$(1+\varepsilon /2)^{-1} < p <-2/\varepsilon .$
}

\bigskip \noindent 
Thangavelu \cite{thangtrans} gave a modification of Kanjin's result by 
replacing the Lebesgue measure $dx$ by $x^{p/4-1/2}dx$ (under the assumption 
$\varepsilon \ge -1/2$). Here we admit more 
general power weights.
\thm{
Let $1<p<\infty,\; \alpha ,\, \beta > -1$ and $ \varepsilon = \min 
\{ \alpha ,\, \beta \} .$  If $\varepsilon \ge 0$,then 
\begin{equation}\label{transplantation}
\| \sum b_k {\cal L}_k^\alpha \| _{L^p_{v(\delta )}}
\le  C \| \sum b_k {\cal L}_k^\beta \| _{L^p_{v(\delta )}}
\end{equation}
for $-1<\delta <p-1 $, where $C$ is a constant independent 
of $f$. If $-1 <\varepsilon <0$, then (\ref{transplantation}) holds for 
$-1-\varepsilon p/2  < \delta < p-1 +\varepsilon  p/2$.
} 

\bigskip \noindent 
Transplantation results were proved for various orthogonal expansions (cf. 
\cite{kan} for a brief exposition).

\medskip \noindent
{\bf Proof.} 
Looking at Kanjin's proof one discerns two lines in the argumentation. \\
The first one consists in pointwise reformulations, estimates, and tools like
the projection formula 
$$L_k^{\mu +\nu }(x) =\frac{\Gamma (k+\mu +\nu +1)}{\Gamma (\nu ) \Gamma (k+ 
\mu +1)} \int _0^1 y^\mu (1-y)^{\nu -1} L^\mu _k (yx) \, dy ,$$
Re$\, \mu >-1,$ Re$\, \nu >0$, discussion of the smoothness properties of the 
involved ``adjusting'' multiplier sequences, verification of the 
Calderon-Zygmund property of a kernel, etc.\\
The second one concerns norm estimates. While Stein's interpolation theorem for 
analytic families of operators and a multiplier criterion of Butzer, Nessel, 
and Trebels \cite{bnt}
do not depend on a particular norm (as long as the hypotheses 
of these theorems are satisfied), there is a dependance of the norm in 
the case of the Calderon-Zygmund theory, Hardy's inequality, and D\l ugosz' 
multiplier theorem. 

\medskip \noindent
If we now follow Kanjin's proof we only have to pay 
attention to the norm estimates and provide the necessary substitutes. We refer 
continuously to the notation used in \cite{kan}. For instance, by $M$ we will 
denote an admissible function, i.e. a positive function $M(\theta ), \, -\infty  
 < \theta < \infty $, that satisfies
$$\sup _{\theta \in {\bf R}} e^{-a|\theta |} \log M(2\theta ) < \infty $$
with some $0<a<\pi $. Also $\varphi (\theta ) =\{ \varphi _n(\theta )\} ,\; 
\theta \in {\bf R},$ will denote the sequence defined by 
$$ \varphi _n(\theta ) =\left( \frac{\Gamma (n+\alpha +1)}{\Gamma (n+\alpha +1 
+i\theta )} \right) ^{1/2} ;$$
the ``adjusting'' operator $T^\beta _{\alpha ,\varphi (\theta )}$ and the 
transplantation operator $T_\alpha ^\beta $ are given by 
$$ T^\beta _{\alpha ,\varphi (\theta )}f \sim \sum 
\varphi _n(\theta ) \langle f,{\cal L}_n^\beta \rangle {\cal L}_n^\alpha , \quad 
\quad T^\beta _\alpha = T^\beta _{\alpha ,\varphi (0)}$$
where $\langle \; ,\; \rangle $ stands for the usual scalar product in 
$L^2({\bf R}_+,dx).$

\medskip \noindent
Let us begin with the case $\alpha , \beta \ge 0 $, temporarily assuming   
 $\max \{ -p/2,-1\} < \delta < \min \{ p/2, p-1 \} $. \\
After some preliminary reductions, similar to those in \cite{kan}, it 
is easy to see that the following proposition (cf. \cite[Proposition 2]{kan}) 
is sufficient to prove Theorem 4.1.
\prop{
Let $\alpha \ge 0,\; 1<p<\infty $, and $k=0,\, 2$. Then 
\begin{equation}\label{reduction}
\| T^{\alpha +k+ i\theta }_{\alpha ,\varphi (\theta )}f \| 
_{L^p_{v(\delta )}}
\le M(\theta ) \| f\| _{L^p_{v(\delta )}}
\end{equation}
for all $\delta $ satisfying
 $\max \{ -p/2,-1\} < \delta < \min \{ p/2, p-1 \} $ with an 
admissible $M$ independently of $f\in C_c^\infty .$
}
To see, for instance, how (\ref{reduction}) implies a weighted analogon of 
\cite[Proposition 1,I]{kan} note that 
$\{ ( \varphi _n(\theta ))^{-1} \} \in wbv_{2,s} ,\; s>\alpha +1,$ whatever 
$\alpha \ge 0$ is and, moreover, 
$$\| \{ ( \varphi _n(\theta ))^{-1}  \} \| _{2,[\alpha +2] } \le C (1+|\theta 
|^{[\alpha +2]}) $$
with $C$ independent of $\theta $. This follows from the calculus for the 
$wbv$-spaces \cite{wbv} since by \cite[Lemma 2]{kan} 
$$ \sup _{x>0} \left| x^j \frac{d^j}{dx^j} \left( \frac{ \Gamma (x+\alpha +1/2 + 
i\theta )}{\Gamma (x+\alpha +1/2)} \right) ^{1/2} \right|
\le C(\alpha ,j)(1+|\theta |^j ) ,\quad j\in {\bf N}_0,\; \alpha >-1/2 ,$$
with $C$ independent of $\theta $. Therefore, by combining (\ref{mspaces}) and 
Theorem 1.1,
$$\| T^{\alpha +k+ i\theta }_\alpha f \| _{L^p_{v(\delta )}}
=\| \sum \langle f,{\cal L}_n^{\alpha +k+i\theta } \rangle {\cal L}_n^\alpha \| 
_{L^p_{v(\delta )}}$$
$$\le C (1+|\theta |^{[\alpha +2]} )\, 
\| T^{\alpha +k+ i\theta }_{\alpha ,\varphi (\theta )}f \| 
_{L^p_{v(\delta )}} \le M(\theta ) \| f\| _{L^p_{v(\delta )}}
$$
provided (\ref{reduction}) holds.

\bigskip \noindent
{\bf Proof of (\ref{reduction}).}
To estimate $\| T^{\alpha + i\theta }_{\alpha ,\varphi (\theta )}f \| 
_{L^p_{v(\delta )}}$ we just follow line by line Section 3 of 
\cite{kan} making use in appropriate places of weighted Hardy's inequality
$$\left( \int _0^\infty \left( \int _x^\infty f(y)\, dy \right) ^p x^\delta dx 
\right) ^{1/p} \le \frac{p}{\delta +1} \left( \int _0^\infty (yf(y))^py^\delta 
dy \right) ^{1/p} $$
valid for $f\ge 0$ and $\delta >-1$ (cf. \cite[p. 272]{sbook}) and of weighted 
inequality for singular integral operators. Recall that the interval $(-1,p-1)$  
characterizes those $\delta $'s for which the function $|x|^\delta $ belongs to 
$A_p({\bf R})$. Therefore we can apply weighted singular integral inequalities.
 
\smallskip \noindent
Exactly the same means are used to estimate 
$\| T^{\alpha + 2+i\theta }_{\alpha ,\varphi (\theta )}f \| 
_{L^p_{v(\delta )}}$, cf. Section 4 of Kanjin's paper, except for 
the fact that we also need a weighted version of \cite[(4.1)]{kan}, i.e., the 
inequality 
\begin{equation}\label{quasi}
\| M^\alpha _\Lambda (f) \| _{L^p_{v(\delta )}}
\le C \| \Lambda \| _{bqc} \| f\| _{L^p_{v(\delta )}}.
\end{equation}
Here one can apply a result of Poiani \cite[Corollary, p. 11]{poi}, which 
gives the uniform boundedness of the Ces\`aro means of order $1$ on weighted 
$L^p(x^\delta dx)$-space, $-1< \delta <p-1$, that is needed in the assumptions
 in 
Butzer, Nessel, and Trebels \cite[Theorem 3.2]{bnt}, and one concludes the 
validity of (\ref{quasi}).

\medskip \noindent
This finishes an outline of necessary changes in Kanjin's proof that have to be 
done in order to prove the weighted transplantation result in the case 
$\alpha , \beta \ge 0 $ and $\max \{ -p/2,-1\} < \delta < \min \{ p/2, p-1 \} $.

\medskip \noindent 
To settle the open cases let us look back at the previous outline and assume 
without loss of generality that $\alpha <\beta $. The 
restrictions on $\alpha , \beta $ as well as on $\delta $ were caused by the 
application of Theorem 1.1. For the proof of (\ref{reduction}) in the case 
$k=2$ and $\theta =0$, 
\begin{equation}\label{red1}
\| T^{\alpha +2}_\alpha f \| 
_{L^p_{v(\delta )}}
\le C \| f\| _{L^p_{v(\delta )}},
\end{equation}
Theorem 1.1 is not needed and a restriction can only come into play by the 
quasi-convexity criterion, i.e. by Poiani's \cite[Corollary, p. 11]{poi} result 
which implies in our case, if $\alpha <0$, that  
$-1-\alpha p/2 < \delta < p-1+\alpha p/2$, hence (\ref{red1}) holds for these 
$\delta$'s.

\medskip \noindent
On the other hand, Theorem 1.1 is true for $-1<\delta =\gamma -\alpha < p-1 $ 
if $\alpha $ is sufficiently large, say $\alpha >A>0$ so that no 
restriction on $\delta $ happens in the transplantation theorem if $\alpha , 
\beta \ge A$. If we choose $N\in {\bf N}_0 $ so large that $\alpha +2N \ge A$, 
the transplantation theorem proved so far, gives 
\begin{equation}\label{reduc2}
\| T_{\alpha +2N}^{\beta +2N} f \| 
_{L^p_{v(\delta )}}
\le C \| f\| _{L^p_{v(\delta )}}.
\end{equation}
The rest of the assertion  now follows by (\ref{red1}),  (\ref{reduc2}), 
duality, the semigroup property of $T_\alpha ^\beta $ (see \cite{kan})
from
$$T_\alpha ^\beta = T_\alpha ^{\alpha +2} \circ 
\cdots \circ T_{\alpha +2N-2}^{\alpha +2N} \circ T_{\alpha +2N}^{\beta +2N}
\circ T_{\beta +2N}^{\beta +2N-2} \circ \cdots \circ T_{\beta +2}^\beta .$$
Hence, Theorem 4.1 is completely established.

\bigskip \noindent

\bigskip \noindent
An immediate consequence of Theorem 4.1 in combination with 
(\ref{mspaces}) yields for Laguerre multipliers  

\coro{ 
Let $\alpha ,\; p, \; \delta $ be as in Theorem 4.1. Then
\begin{equation}\label{transm}
 M^p_{\alpha ,\alpha p/2 +\delta } = {\cal M}^p_{\alpha ,\delta } = 
{\cal M}^p_{0,\delta } =M^p_{0,\delta } . 
\end{equation}
}

\bigskip \noindent
{\bf Remarks.} 1) 
For $\delta =0$ the formula (\ref{transm}) is just Kanjin's \cite{kan} 
multiplier result, 
while for $\delta =p/4 -1/2$ we cover Thangavelu's \cite{thangtrans} statement 
on multipliers.

\medskip \noindent
2) If we choose $\alpha =0$ and $p=2$ then $
M^p_{\alpha ,\alpha p/2} $  and $M^p_{\alpha ,\alpha p/2 +p/4-1/2}$ 
coincide with $M^2_{0,0} =l^\infty $. Since $wbv_{2, s_0} \subset l^\infty ,\; 
s_0>1/2$,
interpolation between Theorem 1.1 for $\alpha =0,\; p>1,\; p\approx 1$ and 
$wbv_{2,s_0} \subset M^2_{0,0}$ along the lines of \cite{cs} 
yields the following improvement of the 
the sufficient multiplier criteria given in  \cite{kan} and 
\cite{thangtrans}: 

\coro{
Let $\alpha >-1$, then
$$wbv_{2,s} \subset M^p_{\alpha ,\alpha p/2} \, , \quad \quad 
wbv_{2,s} \subset M^p_{\alpha ,\alpha p/2 +p/4 -1/2}\, ,\quad  \quad s> \max \{ 
1/p,1-1/p \} ,$$
provided $p$ and $\delta $ satisfy the conditions of Theorem 4.1 where $\delta 
=0$ and $\delta = p/4-1/2$, resp.\\
In particular, if $m(x)$ is a bounded one time  differentiable 
function on ${\bf R}_+$  satisfying
$$ \sup _x |m(x)|^2 +\sup _N \int _N^{2N} |m'(x)|^2 x\, dx \le B^2 ,$$
then, if we set $m_k=m(k)$, there holds
$$ \| T_mf\| _{L^p_{v(\alpha p/2)}} \le CB \| f\| _{L^p_{v(\alpha p/2)}}, \quad 
\quad \| T_mf\| _{L^p_{v(\alpha p/2+p/4 -1/2)}} 
\le CB \| f\| _{L^p_{v(\alpha p/2+p/4-1/2)}} .$$
}

\medskip \noindent 
3) To prove (\ref{lep})
we choose in (\ref{mspaces}) $\gamma = \alpha p(1/p -1/2)$. By the 
transplantation Theorem 4.1 we obtain 
\begin{equation}\label{calm}
 M^p_{\alpha ,\alpha } =  {\cal M}^p_{\alpha ,\alpha p(1/p -1/2)} =
{\cal M}^p_{0,\alpha p(1/p -1/2)} \supset wbv_{2,s_0}, \quad s_0>1 ,
\end{equation}
where the last inclusion follows from Theorem 1.1 provided 
$\max \{ -p/2,-1\} <\alpha p(1/p -1/2)   < \min \{ p/2, p-1 \} $ .
But the latter right inequality leads to a restriction on $p,\; p \le  2$, 
namely 
$$ \alpha p(1/p -1/2)< p-1 \quad \Longleftrightarrow \quad 
\frac{2\alpha +2}{\alpha +2} <p \le 2 ,$$
thus, by duality,  (\ref{lep}) is established.

\medskip \noindent 
4) Interpolation between $wbv_{2,s_0} \subset M^p_{\alpha ,\alpha },\; s_0 >1, 
\alpha \ge 0,$ 
with $p$ as in the preceding remark and 
$l^\infty \subset M^2_{\alpha ,\alpha }$ leads to 
\begin{equation}\label{20}
    wbv_{q,s} \subset M^p_{\alpha ,\alpha } ,\quad s>(2\alpha +2)|1/p-1/2|, 
\quad q<2/s,\quad \frac{2\alpha +2}{\alpha +2} <p <\frac{2\alpha +2}{\alpha } .
\end{equation}
If one applies Corollary 1.2 and the latter criterion upon the multiplier 
sequence $m_{\zeta ,\eta } =\{ m_{\zeta ,\eta }(k) \} ,\; 
m_{\zeta ,\eta }(k) =k^{-\zeta \eta }\exp (ik^\eta ),\; k\neq 0,\; \eta >0,$
one obtains $m_{\zeta ,\eta } \in M^p_{\alpha , \alpha }, \alpha \ge 0,$ 
provided $\zeta > (2\alpha +2)|1/p-1/2|.$

\medskip \noindent
5) Let us conclude with a multiplier criterion for $M^p_{\alpha ,\alpha }$ 
in the case $-1 < \alpha <0$.
Recalling (\ref{calm}),  observing $-1-\alpha p/2 < \alpha -\alpha p/2 <p-1 + 
\alpha p/2$ for all $p, \; 1<p<\infty $, 
and using (\ref{20}) in the case $\alpha =0$, Theorem 4.1 ($\beta =0$) yields 
\coro{
Let $-1<\alpha <0,\; 1<p<\infty $. Then 
$$wbv_{q,s} \subset M^p_{\alpha ,\alpha } \, ,\quad s> 2|1/p-1/2| ,\quad 1/q 
>|1/p-1/2|.
$$
}
\medskip \noindent
This may be looked at as a weak supplement to Corollary 1.2 a) in so far as, if 
one considers the Ces\`aro means of order $\nu $ in the case $\alpha =-1/2$ 
(cf. the discussion to Corollary 1.2), one only obtains  uniform boundedness if 
$\nu >1/4$ for $p>4/3,\; p$ near $4/3$, whereas Muckenhoupt \cite{mu} has shown 
in this instance even the uniform boundedness of the partial sums ($\nu =0$).

\end{document}